\documentclass[12pt]{amsart}
 \usepackage{amsfonts,amssymb,eucal}
\usepackage{amsthm} 
 \usepackage{amsmath} 
\usepackage{amscd}
 \usepackage{latexsym} 
\input{cyracc.def}

\numberwithin{equation}{section}
\usepackage{color}
 \def\ii{\operatorname{i}}
\newcommand{\dd}{\operatorname{d}}
\newcommand{\ad}{{\mbox{\rm ad}}}
\newcommand{\e}{{\mbox{\rm e}}}

\newcommand{\mb}[1]{{\mbox{\boldmath{$#1$}}}}
\newcommand{\mc}[1]{{\mathcal{#1}}}
\newcommand{\got}[1]{{\mathfrak{#1}}}
\newcommand{\db}[1]{{\mathbb{#1}}}

\newcommand{\pa}{\partial}
\newcommand{\R}{\ensuremath{\mathbb{R}}}
\newcommand{\C}{\ensuremath{\mathbb{C}}}
\newcommand{\N}{\ensuremath{\mathbb{N}}}

\newcommand{\Hi}{\ensuremath{\got{H}}}

\newcommand{\Hinf}{\ensuremath{\mathcal{H}^{\infty}}}

\renewcommand{\P}{\ensuremath{\mathbb{P}}}
\newcommand{\Phinf}{\ensuremath{\P (\Hinf )}}
\newtheorem{Remark}{Remark}

\newtheorem{Theorem}{Theorem}

\newtheorem{Proposition}{Proposition}
\newtheorem{lemma}{Lemma}

\theoremstyle{definition}

\newcommand{\oo}{\mathbb{I}}
\newcommand{\un}{{\oo_n}}
\newcommand{\tr}{\ensuremath{{\mbox{\rm{Tr}}}}}%
\begin{document}
\title{Wei-Norman and Berezin's equations of motion on the Siegel-Jacobi disk} 
\author{Stefan  Berceanu}
\address[Stefan  Berceanu]{National
 Institute for Physics and Nuclear Engineering\\
         Department of Theoretical Physics\\
         PO BOX MG-6, Bucharest-Magurele, Romania}
\email{Berceanu@theory.nipne.ro}

\begin{abstract}
We show that the Wei-Norman method applied to describe the evolution  on the
Siegel-Jacobi disk  $\mathcal{D}^J_1=\mathcal{D}_1\times\mathbb{C}^1$, 
where $\mathcal{D}_1$  denotes the Siegel disk,  determined by a hermitian 
Hamiltonian   linear in the generators of the Jacobi group $G^J_1$  and Berezin's
scheme using coherent states give the same equations of  quantum and
classical motion when are expressed in  the coordinates in which the K\"ahler  
two-form $\omega_{\mathcal{D}^J_1} $ can be written as
$\omega_{\mathcal{D}^J_1}=\omega_{\mathcal{D}_1}+\omega_{\mathbb{C}^1}$. 
The Wei-Norman equations on $\mathcal{D}^J_1$  are a particular case
of equations  of motion on the Siegel-Jacobi ball $\mathcal{D}^J_n$
generated by a hermitian Hamiltonian  linear in the generators of the
Jacobi group $G^J_n$ obtained in Berezin's approach based on  coherent 
states on $\mathcal{D}^J_n$.
\end{abstract}

\maketitle
\noindent
\tableofcontents
\newpage
\section{Introduction}\label{intro}

The Jacobi groups -
$G^J_n=H_n\rtimes\text{Sp}(n,\R)_{\C}$ -  where  $H_n$ denotes the 
$(2n+1)$-dimensional Heisenberg group, are unimodular,
nonreductive,  algebraic groups  \cite{ez,bs,TA99, LEE03,gem}  of Harish-Chandra type \cite{satake}. 
The Siegel-Jacobi
domains are nonsymmetric domains associated to the Jacobi groups
by the generalized Harish-Chandra embedding \cite{satake,LEE03, Y02,Y10,gem}.  The Jacobi group  is also an important object in
physics, where  sometimes it is    known  under 
 other names, as 
 {\it Hagen} \cite{hagen},  {\it Schr\"odinger} \cite{ni},    or {\it Weyl-symplectic} group
 \cite{kbw1}. The Jacobi
group describes the   squeezed states  \cite{stol,lu,ho} in Quantum Optics
\cite{mandel,ali,dr}. More references on this subject can be find in \cite{jac1,sbj,sbcg,nou,FC}.

The  Jacobi group  $G^J_n$ has been   studied  in  connection with the
 group-theoretic approach to coherent states 
\cite{perG} in 
\cite{jac1,FC} in the case $n=1$, while the case $n\in\N$ has been
treated  in \cite{mlad,sbj,nou}.  We have attached  to the Jacobi group
$G^J_n$ coherent states  based on
Siegel-Jacobi ball  $\mc{D}^J_n$ \cite{sbj}, which, as set, consists
of the points of 
$\C^n\times\mc{D}_n$.  The non-compact hermitian
symmetric space $ \operatorname{Sp}(n, \R
)_{\C}/\operatorname{U}(n)$ admits a matrix realization  as a bounded
homogeneous domain, the Siegel ball $\mc{D}_n$,
$\mc{D}_n:=\{W\in  M (n, \C ): W=W^t, \un-W\bar{W} > 0\}.$ 
 We have determined the $G^J_n$-invariant
K\"ahler two-form $\omega_{\mc{D}^J_n}$ on $\mc{D}^J_n$  \cite{sbj}, also 
  investigated  by Yang 
\cite{Y10}. In \cite{mlad,sbj,nou}  the $G^J_n$-invariant K\"ahler
two-form $\omega_{\mc{D}_n}(z,W)$,  where $z\in\C^n, W\in \mc{D}_n$, is  written compactly as the sum of two terms, one describing  the homogeneous K\"ahler
two-form   $\omega_{\mc{D}_n}(W)$ on 
$\mc{D}_n$,  the other one is 
$\tr(A^t(\un-\bar{W}W)^{-1}\!\wedge\bar{A})$,  where $ A=\dd z +\dd W\bar{\eta}$, and 
$\eta=(\un-W\bar{W})^{-1}(z+W \bar{z})$,  \cite{mlad,sbj}. We have 
denoted by $FC$ the change of variables $FC: ~
\C^n\times\mc{D}_n$$\ni 
(\eta,W)\rightarrow$$ (z,W)\in\mc{D}^J_n$, $z= \eta-W\bar{\eta}$. We
have shown \cite{nou} that  {\it the $FC$-transform
  is a K\"ahler homogeneous diffeomorphism}, and,  
 when expressed in the variables
$(\eta,W)\in\C^n\times\mc{D}_n$ -- let us call them $FC$-{\it variables} --
the K\"ahler two-form 
$\omega_{\mc{D}_n}(\eta,W)=\omega_{\mc{D}_n}(W)+\omega_{\C^n}({\eta})$. 
We have put in \cite{nou}
this change of variables in connection with the celebrated {\it fundamental
conjecture}  on homogeneous K\"ahler manifolds  of Gindikin and
Vinberg \cite{GV,DN} on the Siegel-Jacobi ball $\mc{D}^J_n$, 
as we did in \cite{FC} for the Siegel-Jacobi disk $\mc{D}^J_1$. Later,
we have underlined in \cite{ber14} that the $FC$-transform has a deep meaning in the
context of Perelomov coherent states  \cite{perG}: it
gives the change of variables from the representation of the normalized
to the un-normalized coherent state vector, as is recalled  in the
Appendix  in \S\ref{app2}.

The equations of motion on the Siegel-Jacobi ball $\mc{D}^J_n$
determined by a hermitian  Hamiltonian   linear in the generators of
the Jacobi group $G^J_n$ were studied in \cite{nou}, generalizing the
results presented  in \cite{jac1,FC},  obtained from the holomorphic differential
representation of the Lie algebra $\got{g}^J_1$  of the Jacobi group.
We recall that linear Hamiltonians in generators of the Jacobi group appear in
  quantum mechanics, as in the case of the  quantum oscillator
acted on by a variable external force \cite{fey,sw,hs} and  in the case of  quantum
dynamics of trapped ions \cite{viorica,ma}.

The aim of this paper is to compare some of the results obtained in
our papers \cite{jac1,FC,nou} concerning the equations of motion  on
the Siegel-Jacobi $\mc{D}^J_n$ determined by a hermitian 
Hamiltonian  linear in the generators of the Jacobi group $G^J_n$ with the
results of the paper \cite{cezar} referring to the equations
of motion on $\mc{D}^J_1$ obtained with the Wei-Norman method
\cite{wei}. In \cite{cezar} the  dynamics on $\mc{D}^J_1$ determined
by a   hermitian  Hamiltonian linear in the generators of $G^J_1$ is considered
in the general framework of Lie systems \cite{lie,lies},  as was  developed
further  in a  geometric approach in \cite{ca1,ca2}. The Wei-Norman
equations for  the  Jacobi group $G^J_1$ in real coordinates have been studied  in
\cite{ca3,ca4}.

In order to establish a correspondence between the formulae of
\cite{cezar} and our notation, we give the following dictionary  (firstly are introduced the symbols used in \cite{cezar}): affine symplectic group $G^{AS}=\text{SL}(2,\R)\rtimes\R^2$ $\leftrightarrow $ Jacobi group
$G^J_1:=H_1\rtimes \text{SU}(1,1)$; extended Poincar\'e disk $\got{M}=\got{D}\times\R^2$ $\leftrightarrow$
Siegel-Jacobi disk $\mc{D}^J_1=\C\times\mc{D}_1$.

Essentially, the  Wei-Norman method (see \cite{wei} and the Appendix
in \S\ref{app1}) consists in representing  the solution
of the equation 
$$\frac{\dd U(t)}{\dd t}= A(t)U(t), \quad U(0)=I, $$
 in the form of product of exponentials
\begin{equation}\label{tsiL}U(t)=\prod_{i=1}^n\exp(\xi_i(t)X_i),
\end{equation}
 where $A$
and $U$ are linear operators, 
\begin{equation}\label{linA}
A(t)=\sum_{i=1}^n\epsilon_i(t)X_i,\end{equation}  
$\epsilon_i(t)$ are scalar functions of
$t$ and $\{X_i\}_{i=1,\dots,n}$ are the generators of a Lie algebra $\got{g}$. 

On the other side, we have developed in \cite{sbcag,sbl}  an algebraic  method to obtain a
representation of a Lie algebra  $\got{g}$ of a Lie group $G$ as first
order holomorphic differential operator  on $M=G/H$  when
$M$ is a hermitian symmetric manifold. Later we have applied    the method
to a larger class of Lie groups, advancing the hypothesis that for the
{\it coherent type groups} \cite{lis,neeb}, i.e. Lie groups for which
the $n$-dimensional homogeneous manifold $M$ admits an holomorphic embedding in a
projective Hilbert space $M\hookrightarrow\Phinf$,  {\it the  generators of
the Lie algebra $\got{g}$ of the Lie group $G$ admit a  
holomorphic differential representation} 
\begin{equation}\label{genPQ}
\got{g}\ni X \mapsto \db{X}(z)=P_X(z) +\sum_iQ^i_X(z)\frac{\pa}{\pa
  z_i}, \end{equation} 
where  $P_X(z)$, and $Q^i_X(z)$ are {\it polynomials} defined on
$M=G/H$. We have verified \cite{jac1,sbj,nou} this hypothesis in the
case of {\it the Jacobi
group} $G^J_n$,  which {\it  is a  coherent type group}
\cite{neeb,ber14}. Following a method advanced in \cite{sbcag,sbl}, 
which uses  Perelomov coherent states \cite{perG} and a {\it
  dequantization} method developed by Berezin \cite{berezin2,berezin1}, we have determined the
equations of motion on $M=\mc{D}^J_n$ when the Hamiltonian $\bf{H}$ is linear
in the generators of the Jacobi group in the case $n=1$ in  
\cite{jac1,FC} and in \cite{nou} for $n\in\N$. In general, for groups $G$  for which the representation \eqref{genPQ}  of the Lie algebra
  $\got{g}$  of the Lie group $G$ is true, the equations of motion on $M=G/H$
depend on the
coefficients $\epsilon_i$ in front of the generators $X_i$ of the group  $G$  which
appear  in  $\mb{H}$ of the form \eqref{linA}  and the
polynomials $Q^i$ which appear in \eqref{genPQ},  as it  is recalled in
Proposition \ref{propvechi}.  To shorten the expression, we call the
equations of motion obtained with this method,  {\it Berezin's
  equations of motion}. We have shown that  for
a Hamiltonian  linear in the generators of $G^J_n$, the motion on
$\mc{D}_n$ is described by a matrix Riccati equation, while the motion
in $z\in\C^n$ is a first order differential equation, with coefficients
depending also on $W\in\mc{D}_n$. It was proved in \cite{FC}  for
$G^J_1$ and in 
\cite{nou} for $G^J_n$, $n\in\N$ that, when the $FC$-transform is applied, the
first order differential equation in the variable $\eta$ becomes
decoupled from the motion on the Siegel ball.   These are {\it  exactly} the
equations of motion obtained in  \cite{cezar}  in the case of the
Siegel-Jacobi disk $\mc{D}^J_1$ using the Wei-Norman method and {\it we
want in the present paper to draw attention to    the fact that
apparently such different methods lead to  the same result}.

The paper is laid out as follows. \S\ref{sec1} recalls  the
definition of the Jacobi algebra $\got{g}^J_1$ adopted  in \cite{jac1}. In \S\ref{repr} the unitary
operators associated with the Jacobi group $G^J_1$ are recalled.  To a 
linear operator  $A$
it is associated the operator  $\hat{\mb{A}}(\xi):=T^{-1}(\xi)
\mb{A}T(\xi)$ \cite{sbcg,sb12,sb13}, where $T(\xi)=D(\alpha)S(w)$, $D(\alpha)$ is the
unitary displacement operator associated to the Heisenberg group
$H_1$,  $S(w)$ is the positive discrete series representation
associated to the group $\text{SU}(1,1)$, and
$\xi=(\alpha,w)\in\C\times\mc{D}_1$.  We take the expressions of 
$\hat{\mb{a}}(\alpha,w) $,  $ \hat{\mb{K}}_{0}(\alpha,w)$ and
$\hat{\mb{K}}_{-}(\alpha,w)$ from  our papers
\cite{sbcg,sb12,sb13}. Then we apply the Wei-Norman method for the
Jacobi group $G^J_1$ in complex, calculating $T^{-1}\frac{\dd T}{ \dd
t}$. In \S\ref{ecMSJ} we determine  the
equations of motion associated to a Hamiltonian  $\mb{H}_0$   linear in the
generators of the Jacobi group $G^J_1$. Following \cite{how}, we
introduce the quasienergy operator $\mb{E}$ associated to the
Hamiltonian $\mb{H}_0$. In \S\ref{Reac} we change the coordinates from
complex to real. The main results of our paper are contained in
Propositions \ref{main} and \ref{conc}. In brief, {\it the Berezin's quantum and classical  equations  of motion on the Siegel-Jacobi disk
  determined by a  hermitian   Hamiltonian linear in the generators of
  Jacobi group  $G^J_1$, expressed in the $FC$-coordinates,  are the same as the equations obtained via the 
  Wei-Norman method. The Wei-Norman equations on $\mc{D}^J_1$  are a
particular case of Berezin's equations of motion on the Siegel-Jacobi ball $\mc{D}^J_n$
generated by a hermitian Hamiltonian  linear in the generators of the Jacobi group
$G^J_n$}. In a  short remark in \S \ref{PHP} are compared the phases
which appears in the method of Wei-Norman \cite{cezar}  and  in Berezin's equations
of motion on the Siegel-Jacobi disk \cite{nou}. For self-containment, in an Appendix in \S\ref{app1}  we
briefly recall the Wei-Norman method. In another Appendix  in \S\ref{app2} are
mentioned the main definitions of coherent states \cite{perG}. In
\S\ref{CLSQ} it is recalled our method for obtaining Berezin's equations of
motion.The construction of coherent states on the Siegel-Jacobi disk
$\mc{D}^J_1$ is summarized in \S \ref{app33} and the equations of
motion on $\mc{D}^J_n$ obtained in \cite{FC,nou} are reproduced in
\S\ref{lasst}, in 
order to make the comparison with the results of \cite{cezar}.
\enlargethispage{1cm}

{\bf Notation}. In this paper the Hilbert space $\Hi$ is endowed with
a scalar product $<\cdot,\cdot>$ 
antilinear in the first argument, i.e. $<\lambda x,y>=\bar{\lambda}<x,y>$,
$x,y\in\Hi,\lambda\in\C\setminus 0$. $\R$, $\C$  and $\N$ denotes the field of real,
complex numbers, respectively the ring of the integers.  We denote the imaginary unit
$\sqrt{-1}$ by $\ii$, and the Real and Imaginary part of a complex
number by $\Re$ and respectively $\Im$, i.e. we have for $z\in\C$,
$z=\Re z+\ii \Im z$, and $\bar{z}=\Re z-\ii \Im z$, but also we use
the notation $cc(z):=\bar{z}$ for  $z\in\C$ or $cc(A)=A^{\dagger}$ for an
operator $A$.  We denote by $M_n(\db{F})$ the set of
$n\times n$ matrices with entries in the field $\db{F}$. If $A\in
M_n(\db{F})$, then 
$A^t$ ($A^{\dagger}$) denotes the transpose (respectively, the
hermitian conjugate) of $A$. $I$  denotes the unit operator,  while
$\un$ denotes the unit matrix  of $M_n(\db{F})$.  If $A\in
M_n(\db{F})$, we denote by $A^s:=\frac{1}{2}(A+A^t)$. If $A$ is a
matrix, then $\tr(A)$ denotes the  trace of the matrix $A$.
We use Einstein convention that repeated indices are implicitly summed.
We
denote the differential by $\dd $. If $\pi$ is an unitary irreducible
representation of a Lie group $G$ with Lie algebra $\got{g}$ on a
complex separable Hilbert space $\Hi$, then we denote for the derived
representation $\mb{X}:=\dd \pi(X)$, $X\in\got{g}$. 
\section{The Lie algebra $\got{g}^J_1$ of the Jacobi group}\label{sec1}

The Heisenberg   group  is the group with the
3-dimensional real  Lie algebra 
 \begin{equation}\label{nr0}\got{h}_1\equiv
<\ii s I+\alpha a^{\dagger}-\bar{\alpha}a>_{s\in\R ,\alpha\in\C} ,\end{equation}
 where   the boson creation  (respectively, annihilation)
operators $a^{\dagger}$ ($a$)  
 verify the canonical commutation relation 
(\ref{baza1}).

We  consider the Lie algebra of the group $\text{SU}(1,1)$:
\begin{equation}\label{nr1}
\got{su}(1,1)=
<2\ii\theta K_0+yK_+-\bar{y}K_->_{\theta\in\R ,y\in\C} , \end{equation} 
where the generators $K_{0,+,-}$ verify the standard commutation relations
(\ref{baza2}).

The Jacobi algebra is defined as the  the semi-direct sum \cite{jac1}
\begin{equation}\label{baza}
\got{g}^J_1:= \got{h}_1\rtimes \got{su}(1,1),
\end{equation}
where $\got{h}_1$ is an  ideal in $\got{g}^J_1$,
determined by the commutation relations \eqref{baza3}, \eqref{baza5}:
\begin{subequations}\label{baza11}
\begin{eqnarray}
& & [a,{{a}^\dagger}]=I\label{baza1}, \\
\label{baza2}
~& & \left[ K_0, K_{\pm}\right]=\pm K_{\pm}~,~ 
\left[ K_-,K_+ \right]=2K_0 , \\
\label{baza3}
& & \left[a,K_+\right]=a^{\dagger}~,~\left[ K_-,a^{\dagger}\right]=a, ~
\left[ K_+,a^{\dagger}\right]=\left[ K_-, a\right]= 0 ,\\
\label{baza5}
& & \left[ K_0  ,~a^{\dagger}\right]=\frac{1}{2}a^{\dagger}, \left[ K_0,a\right]
=-\frac{1}{2}a .
\end{eqnarray}
\end{subequations}
In the conventions of \cite{sbcg},  see equation (3), we have:
\begin{equation}\label{conva}
a=\frac{1}{2\sqrt{\mu}}(P-\ii Q);~ a^{\dagger}=-\frac{1}{2\sqrt{\mu}}(P+\ii Q),  ~[P,Q]=2R.
\end{equation}
which are different of the conventions used in equations (4.14)-(4.16)
in \cite{sbcg}. 
In the convention of  \cite{sbcg}, equation (8), $\mb{P}=\frac{\dd }{\dd x}$, $\mb{Q}=2\ii \mu x$, corresponding to the derived representation of the Heisenberg group, $\mb{R}=\ii \mu \mb{I}$,  $m\in\R$.  The differential realization of \eqref{conva} corresponds to
\begin{equation}\label{difreal}
\mb{a}=\frac{1}{2\sqrt{\mu}}\frac{\dd}{\dd x} + \sqrt{\mu}x:~\mb{a}^{\dagger}=-\frac{1}{2\sqrt{\mu}}\frac{\dd}{\dd x}+ \sqrt{\mu}x.
\end{equation}

\section{Unitary representations associated to the Jacobi group $G^J_1$}\label{repr}
The unitary displacement operator
\begin{subequations}\label{dalpha}
\begin{align}
D(\alpha ) & :=\exp (\alpha {\mb{a}}^{\dagger}-\bar{\alpha}{\mb{a}})\label{da1}\\
~~~& =\exp(-\frac{1}{2}|\alpha
|^2) \exp (\alpha {\mb{a}}^{\dagger})\exp(-\bar{\alpha}{\mb{a}})\label{da2}\\  
~~~& =\exp(\frac{1}{2}|\alpha|^2)  
\exp (-\bar{\alpha} \mb{a})\exp(\alpha\mb{a}^{\dagger})\label{da3}
\end{align}
\end{subequations}
has the composition property
\begin{equation}\label{thetah}
D(\alpha_2)D(\alpha_1)=e^{\ii\theta(\alpha_2,\alpha_1)}
D(\alpha_2+\alpha_1) , 
~\theta(\alpha_2,\alpha_1):=\Im (\alpha_2\bar{\alpha_1}) .
\end{equation}
Note also that
$$D(\alpha)^{\dagger}=D(\alpha)^{-1}= D(-\alpha).$$
Let us denote   by $S$ the unitary irreducible positive discrete
series representation  $D^k_+$  of the group
$\text{SU}(1,1)$ with Casimir operator $C=K^2_0-K^2_1-K^2_2=k(k-1)$,
where $k$ is the Bargmann index for $D^+_k$ \cite{bar47}. 
We introduce the notation $\underline{S}(z)=S(w)$, where 
 $w\in\C,~
|w|<1$ and  $z\in\C\setminus 0$, 
are related by  \eqref{u5}. We have the relations:
\begin{subequations}
\begin{eqnarray}
\underline{S}(z) & := &\exp (z{\mb{K}}_+-\bar{z}{\mb{K}}_-) 
;\label{u1} \\
S(w) & = &  \exp (w{\mb{K}}_+)\exp (\rho
{\mb{K}}_0)\exp(-\bar{w}{\mb{K}}_-)\label{u2}\\
~~~ & =& \exp (-\bar{w}{\mb{K}}_-)\exp (-\rho
{\mb{K}}_0)\exp(w{\mb{K}}_+); \label{u6} \\
w & = &  \frac{z}{|z|}\tanh \,(|z|), ~~\rho =\ln (1-w\bar{w}), z\not=
0,   \label{u5}
\end{eqnarray}
\end{subequations}
and  $w=0$ for $z=0$ in \eqref{u5}.  
Also, it is easy to observe that
$$S(w)^{\dagger}=S(w)^{-1}= S(-w).$$
We  introduce the unitary operator $T(\xi)$:
\begin{equation}\label{txi}
T(\xi)=D(\alpha)S(w), \quad\mc{D}^J_1\ni\xi=(\alpha,w)\in\C\times\mc{D}_1.
\end{equation}
Following \cite{sbcg,sb12,sb13}, for any linear operator $\mb{A}$, 
 we  define the operator
\begin{equation}\label{OpAhat}
\hat{\mb{A}}(\xi):=T^{-1}(\xi) \mb{A}T(\xi), \quad\mc{D}^J_1 \ni\xi=(\alpha,w)\in\C\times\mc{D}_1.
\end{equation}
where $T(\xi)$ was defined in \eqref{txi}.

In  \cite{sb12,sb13}  we have proved that:
\begin{align}\label{kl1}
\hat{\mb{a}}(\alpha,w) & =r (  \mb{a}+w\mb{a}^{\dagger})
+\alpha , \quad r=(1-w\bar{w})^{-\frac{1}{2}}, \\
 \hat{\mb{K}}_{0}(\alpha,w)\! & = r^2[ 
 \bar{w}\mb{K}_{-}\!+\left(  1+|w|^2\right)  \mb{K}_{0}+w\mb{K}_{+}] + \\
 \nonumber & + r \Re[\alpha(
  \mb{a}^{\dagger}+\bar{w}\mb{a})] +\frac{1}{2}|\alpha|^2,\\ 
\label{kl2}
\hat{\mb{K}}_{-}(\alpha,w) & =  r^2[  \mb{K}_{-}+2w\mb{K}_{0}+w^{2}\mb{K}_{+}] + \alpha r (\mb{a}+w\mb{a}^{\dagger}) +\frac{1}{2}\alpha^2.
\end{align}


With formulae \eqref{da2} and \eqref{u2}, we can express \eqref{txi}
as product of exponentials of the generators
as in \eqref{tsiL}, 
where for the Jacobi group $G^J_1$, $n=6$, and the generators
are numbered as 
\begin{equation}\label{x1xn}
\mb{X}_1=\mb{I}; ~\mb{X}_2=\mb{a}^{\dagger};~\mb{X}_3=\mb{a};~ \mb{X}_4=\mb{K}_+;~
\mb{X}_5=\mb{K}_0;~\mb{X}_6=\mb{K}_-,
\end{equation}
while the parameters $\xi_i$  in \eqref{tsiL} are respectively
\begin{equation}\label{xi1n}
\xi_1=-\frac{1}{2}|z|^2;  ~\xi_2=  z;~\xi_3= -\bar{z};  ~\xi_4= w;  ~\xi_5=\ln(1-w\bar{w}) ;  ~\xi_6= -\bar{w} .
\end{equation}
Note that the operator  $T(\xi), ~ \xi=(z,w)\in\mc{D}^j_1$ \eqref{txi}, written as the product \eqref{tsiL} with the generators \eqref{x1xn} and the parameters  \eqref{xi1n},
has the properties expressed in \eqref{tximin}.

Now we apply \eqref{timndot} to the operator  $T(\xi)$ \eqref{txi}  in the
variables $\xi=(z,w)\in\mc{D}^J_1$ expressed with
\eqref{x1xn}, \eqref{xi1n},  taking into account  the commutation
relations \eqref{baza11} of the Lie algebra $\got{g}^J_1$.

For
$$\mb{Y}_2=\e^{-\xi_6\ad\mb{X}_6}e^{-\xi_5\ad\mb{X}_5}e^{-\xi_4\ad\mb{X}_4}e^{-\xi_3\ad\mb{X}_3}\mb{X}_2,$$
we get successively:
\begin{equation}
\begin{split}\label{fY2}
I_1 & = \e^{-\xi_3\mb{X}_3}\mb{X}_2=\mb{a}^{\dagger}-\xi_3;\\
I_2 & = \e^{-\xi_4\mb{X}_4}I_1= I_1; \\
I_3 & = \e^{-\xi_5\mb{X}_5}I_2=
-\xi_3+\e^{-\frac{\xi_5}{2}}\mb{a}^{\dagger};\\
I_4& = \e^{-\xi_6\mb{X}_6}I_3= -\xi_3+\e^{-\frac{\xi_5}{2}}(\mb{a}^{\dagger}-\xi_6\mb{a}).
\end{split}
\end{equation}
For $$\mb{Y}_3=\e^{-\xi_6\ad\mb{X}_6}e^{-\xi_5\ad\mb{X}_5}e^{-\xi_4\ad\mb{X}_4}\mb{X}_3,$$
we get successively:
\begin{equation}
\begin{split}\label{fY3}
e^{-\xi_4\ad\mb{X}_4}\mb{X}_3 & = \mb{a}+\xi_4\mb{a}^
{\dagger};\\
~~e^{-\xi_5\ad\mb{X}_5}\mb{a} & =\e^{\frac{\xi_5}{2}} \mb{a}; \\
~e^{-\xi_5\ad\mb{X}_5}\mb{a}^{\dagger} & = \e^{-\frac{\xi_5}{2}} \mb{a}^{\dagger}.
\end{split}
\end{equation}
For $$\mb{Y}_4=\e^{-\xi_6\ad\mb{X}_6}e^{-\xi_5\ad\mb{X}_5}\mb{X}_4,$$
we get successively:
\begin{equation}
\begin{split}\label{fY4}
J_1 &= e^{-\xi_5\ad\mb{X}_5}\mb{X}_4 = \e^{-\xi_5}
\mb{K}_+;\\
J_2 & = \e^{-\xi_6\ad\mb{X}_6}J_1= \e^{-\xi_5}(
\mb{K}_+-2\xi_6\mb{K}_0+\xi_6^2\mb{K}_-).
\\
\end{split}
\end{equation}
We also have the relations:
\begin{equation}\label{fY5}
\mb{Y}_5=\e^{-\xi_6\ad\mb{X}_6}\mb{X}_5= \mb{K}_0-\xi_6\mb{K}_-.
\end{equation}
Summarizing \eqref{fY2} - \eqref{fY5}, we  obtained the following  expressions
of $\mb{Y}_1 -\ mb{Y}_6$: 
\begin{equation}\label{THEY}
\begin{split}
\mb{Y}_1& =\mb{I};~\mb{Y}_2= -\xi_3+\e^{-\frac{\xi_5}{2}}(\mb{a}^{\dagger}-\xi_6\mb{a});~ \mb{Y}_3=
(\e^{\frac{\xi_5}{2}}-\xi_4\xi_6 \e^{-\frac{\xi_5}{2}})\mb{a}+\xi_4 \e^{-\frac{\xi_5}{2}}\mb{a}^{\dagger}; \\
&  \mb{Y}_4= \e^{-\xi_5}(\mb{K}_++2\bar{\xi}_4\mb{K}_0+\bar{\xi}^2_4\mb{K}_-); ~ \mb{Y}_5= \mb{K}_0+\bar{\xi}_4\mb{K}_-; ~ \mb{Y}_6=\mb{K}_ -.
\end{split}
\end{equation}
Now we introduce the expressions \eqref{THEY} into \eqref{timndot} and we
get  for $T^{-1}\dot{T}$  in the variables $(z,w)\in\mc{D}^J_1$ the expression 
\begin{equation}\label{minTT}
T^{-1}\dot{T}=\ii \Im(z\dot{\bar{z}}) +r[(\dot{z}-\dot{\bar{z}}w)\mb{a}^{\dagger}-cc ] +[\dot{w}r^2\mb{K}_+-cc]
+2\ii \Im (\dot{w}\bar{w})r^2\mb{K}_0 .
\end{equation}

\section{Equations of motion on the Siegel-Jacobi disk
$\mc{D}^J_1$}\label{ecMSJ}
The time-dependent Schr\"odinger equation is expressed as
\begin{equation}\label{SCH}
\mb{H}(t)\psi(t)=\ii\hbar\frac{\dd \psi(t)}{\dd t}.
\end{equation}

As in  \cite{cezar},  we  consider the following family of unitary operators
\begin{equation}\label{uoper}
U(\xi,\varphi):=\exp(-\ii\varphi)T(\xi)
\end{equation}
where $\xi\in\mc{D}^J_1$ and $\varphi$ is a real phase. Let
$\tau=\frac{t}{\hbar}$. In accord with \cite{how},  in \cite{cezar}  it 
was introduced  the quasienergy operator $\mb{E}:=\ii\frac{\dd}{\dd \tau}
-\mb{H}$. With \eqref{OpAhat},  we get:
\begin{equation}\label{Ehat}
\hat{\mb{E}}(\xi,\varphi)=\frac{\dd\varphi }{\dd \tau}\mb{I}+\ii T(\xi)^{-1}\dot{T}(\xi)-\hat{\mb{H}}(\xi).
\end{equation}

In the notation of \cite{jac1,FC}, we consider a hermitian Hamiltonian linear in the generators of the Jacobi
group $G^J_1$:
\begin{equation}\label{guru}
\mb{H}_0 = \epsilon_a\mb{a} +\bar{\epsilon}_a\mb{a}^{\dagger}
 +\epsilon_0 {\mb{K}}_0 +\epsilon_+{\mb{K}}_++\epsilon_-{\mb{K}}_-  ,~~
\bar{\epsilon}_+=\epsilon_-, ~ \epsilon_0=\bar{\epsilon_0}.
\end{equation}
With equations \eqref{kl1}-\eqref{kl2}, we calculate
$\hat{\mb{H}}_0(\xi)$, where  $\xi=(\alpha,w)\in\C\times\mc{D}_1$: 
\begin{equation}\label{guruhat}
\hat{\mb{H}}_0 (\alpha,w) =I_0+C_1\mb{a}^{\dagger}+\bar{C}_1\mb{a}+C_0\mb{K}_0+C_+\mb{K}_++ \bar{C}_+\mb{K}_-,
\end{equation}
where the coefficients in  \eqref{guruhat} have the values:
\begin{equation}
\begin{split}\label{ioc}
I_0 & =\epsilon_a\alpha+\bar{\epsilon_a}\bar{\alpha}+\frac{1}{2}(\epsilon_0|\alpha|^2+\epsilon_-\alpha^2+
\epsilon_+\bar{\alpha}^2),\\
\frac{C_1}{r}  & = \bar{\epsilon}_a+\epsilon_{a}w+\frac{\epsilon_0}{2}(\alpha+\bar{\alpha}w)
+\epsilon_-\alpha w+\epsilon_+\bar{\alpha},\\
\frac{C_0}{r^2} & =\epsilon_0(1+|w|^2)+2(\epsilon_-w+\epsilon_+\bar{w}),\\
\frac{C_+}{r^2} & =\epsilon_0w+\epsilon_-w^2+\epsilon_+. \\
\end{split}
\end{equation}
\section{The Jacobi group $G^J_1$ in real coordinates}\label{Reac}

In \cite{cezar} it is used as basis of the Lie algebra $\got{g}^J_1$ 
the real basis from \cite{bs,sbcg}.

We now list the relations between the operators used in the paper
\cite{cezar}  and the generators \eqref{baza11} of the Lie algebra $\got{g}^J_1$: 
\begin{equation}\label{n1n2}\mb{N}_1=\mb{a}+\mb{a}^{\dagger};\quad\mb{N}_2=\ii (\mb{a}-\mb{a}^{\dagger}),
\end{equation}
with  the inverse
\begin{equation}
\mb{a}=\frac{1}{2}(\mb{N}_1-\ii \mb{N}_2); \quad \mb{a}^{\dagger}=\frac{1}{2}(\mb{N}_1+\ii \mb{N}_2);
\end{equation}
\begin{equation}\label{k1k2}
\mb{K}_1=\frac{1}{2}(\mb{K}_++\mb{K}_-); \quad \mb{K}_2=\frac{1}{2\ii}(\mb{K}_+-\mb{K}_-),
\end{equation}
and the inverse
\begin{equation}
\mb{K}_+=\mb{K}_1+\ii\mb{K}_2;\quad \mb{K}_-=\mb{K}_1-\ii\mb{K}_2.
\end{equation}
In \cite{cezar} it was considered the unitary operator
\begin{equation}\label{TCEZ}
T(\xi)=D(x,y)S(u,v),\quad \xi=(u,v,x,y)\in\got{M},
\end{equation}
where
\begin{equation}\label{DScezar}
D(x,y)=\exp(\ii y\mb{N}_1+\ii x\mb{N}_2), \quad S(u,v)=\exp(\ii k_1\mb{K}_1 +\ii k_2\mb{K}_2),
\end{equation}
\begin{equation}\label{KK12}
k_1=\frac{v}{2s}\ln\frac{1+s}{1-s}, \quad k_2=\frac{u}{2s}\ln\frac{1+s}{1-s},\quad s=(u^2+v^2)^{\frac{1}{2}} .
\end{equation}
The correspondence between the real parameters in the representation \eqref{DScezar} and the complex parametrization \eqref{txi} is given by the relations
\begin{equation}\label{corT}
\alpha= x+\ii y; \quad w=u+\ii v. 
\end{equation}
In \cite{cezar}  the Hamiltonian \eqref{guru} was written down as:
\begin{equation}\label{hcezar}
\mb{H}_0=2\varepsilon_0\mb{K}_0+2\varepsilon_1\mb{K}_1+ 2\varepsilon_2\mb{K}_2 + 2\nu_1\mb{N}_1  + 2\nu_2\mb{N}_2. 
\end{equation}
 The correspondence of the real and complex coefficients of the
Hamiltonians \eqref{hcezar}  and \eqref{guru} is (see also Proposition \ref{EQLIN1})
\begin{equation}\label{corc}
\epsilon_a=\nu_1+\ii \nu_2; ~\epsilon_0=2 \varepsilon_0; ~ \epsilon_+=\varepsilon_1-\ii \varepsilon_2, 
\end{equation}
\begin{equation}\label{corcr}
\nu_1=a; ~\nu_2=b; ~\varepsilon_0=p;~ \varepsilon_1=m;~\varepsilon_2=n. 
\end{equation}

We express  $T^{-1}\dot{T}$ \eqref{minTT}  given in the complex
coordinates $(z,w)\in\mc{D}^J_1$ in the real coordinates $(x,y;u,w)$,
where 
$z=x+\ii y$, $w=u+\ii v$, and in the real operators \eqref{n1n2}, \eqref{k1k2}, and we get:
\begin{equation}\label{minTT1}
\begin{split}
-\ii T^{-1}\frac{\dd T}{\dd \tau } & =(x\dot{y}-y\dot{x})\mb{I}+r[(1+u)\dot{y}-\dot{x}v]\mb{N}_1+r[(1-u)\dot{x}-\dot{y}v]\mb{N}_2\\
~~~&  +2r^2[(\dot{v}u-\dot{u}v)\mb{K}_0+\dot{v}\mb{K}_1+\dot{u}\mb{K}_2].
\end{split}
\end{equation}
Now we express  the operator $\hat{\mb{H}}_0(\alpha,w)$ \eqref{guruhat} in the
operators \eqref{n1n2}, \eqref{k1k2} in real coordinates $(x,y;u,v)$, where
$\alpha=x+\ii y$, $w=u+\ii v$, and we get:
\begin{equation}\label{hatR}
\hat{\mb{H}}_0(x,y;u,v)=D_0+D_1\mb{N}_1+ D_2\mb{N}_2+F_0\mb{K}_0+ F_1\mb{K}_1+F_2\mb{K}_2.
\end{equation}
We find the following values of the coefficients appearing in  \eqref{hatR}:
\begin{equation}
\begin{split}\label{h2r}
D_0 &= 2(\nu_1x-\nu_2y)+\frac{\epsilon_0}{2}(x^2+y^2)+\varepsilon_1(x^2-y^2)-2\varepsilon_2 xy,\\
\frac{D_1}{r} &=\nu_1(1+u)-\nu_2v+\frac{\epsilon_0}{2}[x(1+u)+yv]\\~~~&+\varepsilon_1[x(1+u)-yv] -\varepsilon_2[xv+(1+u)y],
\\-\frac{D_2}{r} &= \nu_1v+\nu_2(u-1)+\frac{\epsilon_0}{2}[xv+y(1-u)] \\~~~&\varepsilon_1[xv+y(u-1)]
+\varepsilon_2[x(u-1)-yv],\\
\frac{F_0}{r^2} & =\epsilon_0(1+u^2+v^2)+4(\varepsilon_1u-\varepsilon_2v),\\
\frac{F_1}{2r^2} &
=\epsilon_0u+\varepsilon_1(u^2-v^2+1)-2\varepsilon_2uv,\\
-\frac{F_2}{2r^2} & =\epsilon_0v+2\varepsilon_1uv+\varepsilon_2(-1+u^2-v^2).
\end{split}
\end{equation}
We introduce \eqref{hatR} and \eqref{minTT1} into \eqref{Ehat} and we
get the expression
\begin{equation}\label{Ehat1}
\hat{\mb{E}} (x,y;u,v)= G_0\mb{I}+G_1\mb{N}_1+G_2\mb{N}_2+H_0\mb{K}_0+H_1\mb{K}_1+H_2\mb{K}_2,
\end{equation}
where:
\begin{equation}
\begin{split}\label{G0H}
G_0 & = \dot{\varphi} +y\dot{x}-x\dot{y}-2(\nu_1x-\nu_2y)-\varepsilon_0(x^2+y^2) -\varepsilon_1(x^2-y^2)+2\varepsilon_2xy,\\
-\frac{G_1}{r} &= (1+u)\dot{y}-\dot{x}v+\nu_1(1+u)-\nu_2v+\varepsilon_0[x(1+u)+yv]\\ ~~~&+\varepsilon_1[x(1+u)-yv]
-\epsilon_2[y(1+u)+xv],\\
\frac{G_2}{r} &= -(1-u)\dot{x}+\dot{y}v+\nu_1v+\nu_2(u-1)+\varepsilon_0[y(1-u)+xv] \\ ~~~&+\varepsilon_1[xv+y(u-1)] +\varepsilon_2[x(u-1)-yv],\\
-\frac{H_0}{2r^2} & = \dot{v}u-\dot{u}v+\varepsilon_0(1+u^2+v^2)+2(\varepsilon_1u-\varepsilon_2v),\\
-\frac{H_1}{2r^2} & =  \dot{v}+\varepsilon_0u+\varepsilon_1(u^2-v^2+1)-2\varepsilon_2uv,\\
\frac{H_2}{2r^2} & = -\dot{u}+\varepsilon_0v+2\varepsilon_1uv+\varepsilon_2(u^2-v^2-1).
\end{split}
\end{equation}
Identifying the coefficients of $\mb{N}_1$, $\mb{N}_2$, and
respectively  $\mb{K}_1$, $\mb{K}_2$, we get the equations of motion
in real coordinates
\begin{equation}\label{xyec}
\begin{split}
\dot{x} & = -\varepsilon_2x+(\varepsilon_0-\varepsilon_1)y-\nu_2,\\
\dot{y} & = -(\varepsilon_0+\varepsilon_1)x+\varepsilon_2y-\nu_1; 
\end{split}
\end{equation}
\begin{equation}\label{uvec}
\begin{split}
\dot{u} & = 2v(\varepsilon_1u+\varepsilon_0)-\varepsilon_2(1-u^2+v^2),\\
\dot{v}&  =  2u(\varepsilon_2v-\varepsilon_0)-\varepsilon_1(1+u^2-v^2). 
\end{split}
\end{equation}
Introducing \eqref{uvec} into the equation of $H_0$ in \eqref{G0H}, we have
\begin{equation}\label{H00}
-\frac{H_0}{2}=\varepsilon_0+\varepsilon_1u-\varepsilon_2v.
\end{equation}
Introducing \eqref{xyec} into expression of $G_0$ in  \eqref{G0H}, we have
\begin{equation}\label{GOOH}
G_0=\dot{\varphi}-(\nu_1x-\nu_2y).
\end{equation}
Starting from the complex representation, we have regained  the
results of Section 3 in \cite{cezar} and also our results from
\cite{nou}, reproduced in Proposition \ref{EQLIN1}:
\begin{Proposition}\label{main}
If we consider $\Phi\in\Hi$ such that $\mb{K}_0\Phi=k\Phi$, then
\begin{equation}\label{SCHCEZ}
\Psi(\xi,\varphi)=U(\xi,\varphi)\Phi= \e^{-\ii\phi}T(\xi)\Phi\end{equation} is a solution of the time
dependent Schr\"odinger equation \eqref{SCH} corresponding to the
hermitian Hamiltonian \eqref{guru} (or \eqref{hcezar}) on the
Siegel-Jacobi disk $\mc{D}^J_1$, where $x,y\in\R$ verify \eqref{xyec},
$(u,v)\in\mc{D}_1$ verify \eqref{uvec}, while the phase $\varphi$ in \eqref{uoper}
verifies 
\begin{equation}\label{phicez}
\dot{\varphi}=\nu_1x-\nu_2y+2k(\varepsilon_0+\varepsilon_1u-\varepsilon_2v).
\end{equation}

The motion \eqref{uvec} on the Siegel disc $\mc{D}_1$
in the 
  complex variable $w=u+\ii v$ is described by the Riccati equation 
\begin{equation}\label{wcom}
\ii \dot{w}= \epsilon_++\epsilon_0 w+\epsilon_-w^2.
\end{equation}

The equations of motion \eqref{xyec} in 
 the complex variable $\eta= x+\ii y$ reads
\begin{equation}\label{etaec}
\ii \dot{\eta}=\bar{\epsilon}_a+\epsilon_+\bar{\eta}+\frac{\epsilon_0}{2}\eta,
\end{equation}
\end{Proposition}

We remark that 
 the Riccati equation \eqref{wcom} on $\mc{D}_1$ obtained with the Wei-Norman method
 coincides with the Berezin's equation of motion {\rm{(4.8b)}}  or
 {\rm{(4.10b)}}   in \cite{FC}, with the
difference of notation $\epsilon_+\leftrightarrow
\epsilon_-=\bar{\epsilon}_+$.  
The equation \eqref{etaec} for $\eta\in\C$ obtained with the Wei-Norman method is the
Berezin's equation of motion {\rm{(4.10a)}}    in
\cite{FC},  with the
correspondence $\epsilon_a\leftrightarrow
\bar{\epsilon}_a$,  $\epsilon_+\leftrightarrow
\epsilon_-=\bar{\epsilon}_+$, see Propositions \ref{EQLIN1} and  \ref{POYT}.

We have proved 
\begin{Proposition}\label{conc}The quantum and classical    Berezin's equations
  of motion on the Siegel-Jacobi disk
  determined by a  hermitian   Hamiltonian,  linear in the generators of
  Jacobi group  $G^J_1$,  expressed in the $FC$-coordinates,  are the same
  as the equations obtained applying the Wei-Norman method.

The equations of motion \eqref{wcom} on $\mc{D}_1$,
  \eqref{etaec} on $\C$  and \eqref{xyec} on $\R^2$,  determined by
 the  hermitian  Hamiltonian \eqref{guru} or \eqref{hcezar}  linear in  the generators of the Jacobi group $G^J_1$ obtained with the
  Wei-Norman method are a particular case of the Berezin's equations of motion
  \eqref{Mhip2} on $\mc{D}_n$, \eqref{hipPRT1} on $\C^n$ and
  respectively \eqref{LINe} on $\R^{2n}$,  determined by the 
  hermitian Hamiltonian \eqref{HACA} linear in the generators of the
  Jacobi group $G^J_n$.

\end{Proposition}

\section{Phases}\label{PHP}
In the paper \cite{cezar} it was calculated the phase $\varphi$ which
appears  in the solution \eqref{uoper} for which we have find the
equation \eqref{phicez}. In the Berezin's approach to the equations of
motion, the solution of the time-dependent Schr\"odinger equation
\eqref{SCH}
differs from the solution parametrized by Perelomov coherent states by
a phase $\phi$ as it is recalled in  Proposition \ref{propvechi}
proven in \cite{sbcag,FC,nou}. 
Comparing the Wei-Norman solution of the Schr\"odinger equation
\eqref{SCH} for the Hamiltonian \eqref{guru}  (or \eqref{hcezar}) with
the solution \eqref{slSCH} under the form \eqref{legPH}, we see that

\begin{Remark}\label{REMPS}
The phases $\varphi$ used in the Wei-Norman method  associated to the
quasi-energy  operator \eqref{Ehat} and the phase $\phi$ which
appears in Berezin's approach are different: 
\begin{equation}\label{IMPDIF}
-\varphi(\xi)=\phi(\xi)-\frac{1}{2}\Im(\omega\bar{\alpha}^2),\quad \xi=(\alpha,w)\in\mc{D}^J_1.
\end{equation}
If we introduce the equations of motions  on $\mc{D}^J_1$
into the-$\tau$ derivative of \eqref{IMPDIF},  we get
\begin{equation}\label{fiPHI}
-\dot{\phi}=\dot{\varphi}
-\nu_1(ux+vy)+\nu_2(uy-vx)+\frac{1}{4}(\epsilon_-\bar{z}^2+\epsilon_+z^2), \end{equation} 
where $z=\alpha-w\bar{\alpha}$ and $(\alpha,w)$ are the coordinates \eqref{txi} on
$\mc{D}^J_1$. This is exactly the formula obtained for $\dot{\phi}$
summing up  the explicit
expressions of the dynamical and Berry phases  \eqref{realHH} and \eqref{FFV}.
\end{Remark}
We verify the last part of the Remark. 
If we add the  expression of the Berry phases given by
\eqref{FFV} in which we introduce the equations of motion
\eqref{qqqN1} determined by the Hamiltonian  \eqref{guru} and dynamic
phase \eqref{realHH}, we get  the value of  the
phase $\phi$  in the complex  variables
$(z,w)\in\C\times\mc{D}_1$, $z=\alpha-w\bar{\alpha}$:
\begin{subequations}\label{pp}
\begin{align}-{\phi} & =k{\phi}_1 +{\phi}_0,\\
\dot{\phi}_1 & = \epsilon_0+\epsilon_-\bar{w}+\epsilon_+w,\label{pp1}\\
\dot{\phi}_0 & = \frac{1}{4}(\epsilon_-\bar{z}^2+\epsilon_+z^2)+
\frac{1}{2} (\epsilon_a\bar{z}+\bar{\epsilon}_az) \label{pp0}. 
\end{align}
\end{subequations}
When we express $z$ in the real and imaginary part in the coordinates \eqref{corT}
$(x,y,u,v)$, we have for $z$ appearing in \eqref{pp0} the value
$$z=(1-u)x-yv+\ii [(1+u)y-vx],$$
while $$\dot{\phi}_1=2(\varepsilon_0+\varepsilon_1u +\varepsilon_2 v),$$
i.e. the expression  multiplying   $k$ in    \eqref{phicez},
with the correspondence $\epsilon_+\leftrightarrow
\epsilon_-=\bar{\epsilon}_+$, $\epsilon_a\leftrightarrow
\bar{\epsilon}_a$.  \\[3ex]

{\bf In conclusion}, in \cite{nou,FC,ber14} we have underlined the utility 
of the $FC$-coordinates,  the geometric  significance of the 
$FC$-transform in the context of the fundamental conjecture and also
its relevance for coherent states on the Siegel-Jacobi ball. In the
present  paper {\it we have shown that Berezin's equations of motion
on the Siegel-Jacobi disk   expressed  in the $FC$-coordinates, 
generated by a linear Hamiltonian in the generators of th Jacobi group 
are identical with the equations of
motion furnished by the Wei-Norman method}. All the calculation in the
present paper refers to the motion on the Siegel-Jacobi disk,  generated
by  a linear Hamiltonian in the generators of the Jacobi group $G^J_1$, but we
believe  that  Proposition  \ref{conc} is  true in more general
situations,  for some Lie groups which are semidirect product. We have also underlined that the phases which appears
in the two  methods are different.  Remark
\ref{REMPS}  is also  a direct  check of the correctness of our calculation
in \cite{FC} and compatibility with the  calculation in \cite{cezar}.

\section{Appendix: The Wei-Norman method}\label{app1}

We use the following convention of notation for noncommuting operators:
\begin{equation}\label{conv}
\prod_{i=1}^nA_i:=A_1\dots A_n; \quad\prod_{i=n}^1A_i:=A_n\dots A_1.
\end{equation}

Let $\xi=(\xi_1,\dots,\xi_n)$ be  some parameters. We consider an unitary
operator which can be expressed in the basis $\{X_i\}_{i=1,\dots,n}$
of the Lie algebra $\got{g}$ as in \eqref{tsiL}. 
Then we have
\begin{equation}\label{tximin}
U^{-1}(\xi)=U^{\dagger}(\xi)=U(-\xi)=\prod_{j=n}^1\exp(-\xi_j\mb{X}_j).
\end{equation}
Let $X,Y$ be the free generators of a ring $R$. The Baker-Hausdorff
formula (see  \cite{mag} for a proof)  reads :
\begin{equation}\label{BCH}\begin{split}
e^XYe^{-X}&=e^{\ad X}Y=\sum_{n=0}^{\infty}\frac{(\ad X)^n}{n!}Y\\
~~~&=Y+[X,Y]+\frac{1}{2!}[X,[X,Y]]+\dots
+\frac{1}{n!}[\underbrace{[X,[X,\dots,[}_{n\text{~brackets}}X,Y\underbrace{]\dots]}_{n}+\dots  .
\end{split}
\end{equation}
Now let us consider that the parameters $\xi$ depend on a variable,
let call it $t$. For any $t$-dependent operator $A$, we  denote
$\dot{A}=\frac{\dd A}{\dd t}$. Then we can calculate the derivative of \eqref{tsiL} as
\begin{equation}\label{dotT}
\dot{U}(\xi)=\sum_{j=1}^n \dot{\xi}_j\left(\prod_{k=1}^{j-1}[\exp(\xi_k\mb{X}_k)]\mb{X}_j\prod_{l=j}^n[\exp(\xi_l\mb{X}_l)]\right).
\end{equation}
Let us consider a linear operator $A(t)$ of the form \eqref{linA}, 
where $a_i(t)$ are scalar functions of $t$.
It can be   proved  \cite{wei} that:
\begin{lemma}If $\{\mb{X}_1,\dots,\mb{X}_n \}$ is  a basis of a Lie algebra
  $\got{g}$, then
\begin{equation}\label{lemW}
\left(\prod_{j=1}^r\exp(\xi_j\mb{X}_j)\right)\mb{X}_i\left(\prod_{j=r}^1\exp(-\xi_j\mb{X}_j)\right)=\sum_{k=1}^n\eta_{ki}\mb{X}_k,~i,
r=1,2,\dots,n,
\end{equation}
where $\eta_{ki}=\eta_{ki}(\xi_1,\dots,\xi_r)$ is an analytic
function of $\xi_1,\dots,\xi_r$. 
\end{lemma}
The Wei-Norman method  \cite{wei} is expressed in the theorem:
\begin{Theorem}
If $A(t)$ is  given by \eqref{linA},  then there exists a neighborhood
of $t=0$ in which the solution of the equation
\begin{equation}\label{dotU}
\dot{U}(t)=A(t)U(t), \quad U(0)=I
\end{equation}
may be expressed in the form \eqref{tsiL}, where the $\xi_i(t)$ are
scalar functions of $t$. Moreover, $\xi^t:=(\xi_1,\dots \xi_n)$ satisfy the first order 
differential equation
\begin{equation}\label{solwei}
\eta\dot{\xi}=\epsilon,
\end{equation}
 which depend only on the Lie algebra $\got{g}$
and the $\epsilon(t)$'s. $\eta=(\eta_{ki})_{k,i=1,\dots,n}$ is the matrix of coefficients of \eqref{linA}, while $\epsilon$
denotes the vector with coefficients which appear in \eqref{linA},  $ \epsilon^t=(\epsilon_1,\dots,\epsilon_n)$. 
\end{Theorem}
In \cite{wei} it was proved that the representation \eqref{tsiL} is
global for any solvable Lie algebra $\got{g}$ and for any $2\times 2$
system of equations.

In our calculation  in \eqref{Ehat}, instead of $\dot{T}T^{-1}$, we need
$T^{-1}\dot{T}$. 
With the convention \eqref{conv} and the Baker-Hausdorff formula \eqref{BCH}, we obtain
for $U$ defined in \eqref{tsiL}:
\begin{equation}\label{timndot}
U^{-1}(\xi)\dot{U}(\xi)=\sum_{i=1}^n\dot{\xi}_i\mb{Y}_i, \quad \mb{Y}_i= \left(\prod_{k=n}^{i+1}[\exp(-\xi_k\ad \mb{X}_k)]\right)\mb{X}_i.
\end{equation}

\section{Appendix: Berezin's equations on motion}\label{app2}
We  consider the triplet $(G, \pi , \got{H} )$, where $\pi$ is
 a continuous, unitary, irreducible 
representation
 of the  Lie group $G$
 on the   separable  complex  Hilbert space \Hi~   \cite{perG}.

We
 introduce the normalized (unnormalized) vectors  $\underline{e}_x$
(respectively, $e_z$) defined on $G/H$
\begin{equation}\label{bigch}
\underline{e}_x=\exp(\sum_{\phi\in\Delta_+}x_{\phi}\mb{X}_{\phi}^+-\bar{x}_{\phi}\mb{X}_{\phi}^-)
e_0, ~
e_z=\exp(\sum_{\phi\in\Delta_+}z_{\phi}\mb{X}_{\phi}^+)e_0, 
\end{equation}
where
$e_0$ is the extremal weight vector of the representation $\pi$, $\Delta_+$ are the positive roots of the Lie algebra $\got{g}$ of $G$,
and $X_\phi,\phi\in\Delta$,  are the  generators. 
$\mb{X}^+_{\phi}$ ($\mb{X}^-_{\phi}$)
corresponds to the positive (respectively, negative) generators. See
details in  \cite{perG,sb6}.

We denote by $FC$ the change of variables
$x\rightarrow z$ in formula \eqref{bigch} such that
\begin{equation}\label{etild}
\underline{e}_{x}=\tilde{e}_z, ~    \tilde{e}_z  :=
  (e_z,e_z)^{-\frac{1}{2}}e_z, ~z=FC(x). 
\end{equation}
The reason for calling the transform \eqref{etild}  a $FC$-transform  {\it
  (fundamental conjecture)}  \cite{GV,DN} is explained in Proposition
3 in \cite{ber14}. For a concrete example of $FC$-transform, see \eqref{csv}.

 \subsection{Berezin's approach to classical motion and quantum evolution}\label{CLSQ}

Let $M=G/H$ be a homogeneous  manifold with a $G$-invariant K\"ahler two-form
$\omega$
\begin{equation}\label{kall}
\omega(z)=\ii\sum_{\alpha\in\Delta_+} g_{\alpha,\beta}  d z_{\alpha}\wedge
d\bar{z}_{\beta}, ~g_{\alpha,\beta}=\frac{\pa^2}{\pa
  z_{\alpha} \pa\bar{z}_{\beta}} \ln  <e_z,e_z>. 
\end{equation}

Passing on from the dynamical system problem
 in the Hilbert space $\Hi$ to the corresponding one on $M$ is called
sometimes {\it dequantization}, and the dynamical system on $M$ is a classical
one \cite{sbcag,sbl}. Following Berezin \cite{berezin2,berezin1}, the
motion on the classical phase space can be described by the local
equations of motion
$\dot{z}_{\alpha}=\ii \left\{\mc{H},z_{\alpha}\right\},
  ~\alpha \in \Delta_+ $, where $\mc{H}$ is
  the classical Hamiltonian  ({\it the covariant 
  symbol})
\begin{equation}\label{clH}\mc{H}=<e_z,e_z>^{-1}<e_z|\mb{H}|\e_z>\end{equation}
  attached to
  the quantum Hamiltonian $\mb{H}$, and the Poisson bracket is
  introduced using the matrix $g^{-1}$. 

We consider an algebraic Hamiltonian linear in the generators
${\mb{X}}_{\lambda} $ of the
group of symmetry $G$
\begin{equation}\label{lllu}
\mb{H}=\sum_{\lambda\in\Delta}\epsilon_{\lambda}{\mb{X}}_{\lambda} .
\end{equation}
We look for the solution of the Schr\"odinger equation of motion
\eqref{SCH} generated by the Hamiltonian \eqref{lllu} 
as \begin{equation}\label{slSCH}
\psi(t)=\e^{\ii \phi}\tilde{e}_z,
\end{equation}
where $\tilde{e}_z$ is the normalized Perelomov coherent state vector
defined in \eqref{bigch}, \eqref{etild}.

We extract from  \cite{sbcag,sbl,nou}  the following Proposition: 
\begin{Proposition}\label{propvechi}
The  classical motion and quantum evolution  generated
by the linear hermitian  Hamiltonian  \eqref{lllu}  are described
by Berezin's   equations
of motion  on $M=G/H$
\begin{equation}\label{moveM}
{\ii\dot{z}_{\alpha}=\sum_{\lambda\in\Delta}\epsilon_{\lambda}Q_{\lambda
,\alpha}},~\alpha\in\Delta_+ , 
\end{equation}
where the differential action corresponding to the operator
$\mb{X}_{\lambda}$ in \eqref{lllu} can be expressed in a local
system of coordinates as a holomorphic  first order differential
operator with polynomial coefficients  ($\pa_{\beta}=\frac{\pa}{\pa
  z_{\beta}}$),
\begin{equation}\label{VBC}\db{X}_{\lambda}=P_{\lambda}+\sum_{\beta\in\Delta_+}Q_{\lambda,
    \beta}\partial_{\beta}, \lambda\in\Delta.
\end{equation}
The phase $\phi$ in \eqref{slSCH} is given by the sum
$\phi=\phi_D+\phi_B$
of the dynamical and Berry phases, where
\begin{subequations}
\begin{align}
\dot{\phi}_D & = -\mc{H}(t); \\
\dot{\phi}_B & = \frac{\ii}{2}\sum_{\alpha\in\Delta_+}(\dot{z}_{\alpha}\pa_{\alpha}-
 \dot{\bar{z}}_{\alpha}\bar{\pa}_{\alpha})\ln <e_z,e_z>. 
\end{align}
\end{subequations}
\end{Proposition}

\subsection{Coherent states on the  Siegel-Jacobi disk $\mc{D}^J_1$}\label{app33}

 We impose to the cyclic vector $e_0$ to verify simultaneously
 the conditions \cite{jac1}
\begin{equation}\label{cond}
\mb{a}e_0  =  0, ~
 {\mb{K}}_-e_0  =  0,~
{\mb{K}}_0e_0  =  k e_0;~ k>0, 2k=2,3,... , 
\end{equation}
and we have considered in  the last relation in (\ref{cond}) the positive  discrete series
representations $D^+_k$ of $\text{SU}(1,1)$  \cite{bar47}.

 Perelomov's coherent state   vectors   associated to the group $G^J_{1}$ with 
Lie algebra the Jacobi algebra $\got{g}^J_{1}$, based on Siegel-Jacobi
disk $ \mathfrak{D}^J_{1}  =  H_1/\R\times \text{SU}(1,1)/\text{U}(1)$
$ =\C\times\mc{D}_1$,  
are defined as 
\begin{equation}\label{csu}
e_{z,w}:=\e^{z\mb{a}^{\dagger}+w{\mb{K}}_+}e_0, ~z,w\in\C,~ |w|<1 .
\end{equation}
We introduce also the normalized    ({\it squeezed})     CS-vector
(see also \cite{stol}) 
\begin{equation}\label{sqz}
\underline{e}_{\xi}:= T(\xi)  e_0=D(\alpha )S(w) e_0, \quad \xi=(\alpha, w)\in\C\times\mc{D}_1.\end{equation}
The normalized squeezed state vector   and the
un-normalized 
 Perelomov's coherent state vector
are related by the relation (see \cite{jac1})
\begin{equation}\label{csv}
\underline{e}_{\eta, w} = (1-w\bar{w})^k
\exp (-\frac{\bar{\eta}}{2}z)e_{z,w},~ z=\eta-w\bar{\eta}.
\end{equation}
We recall \cite{jac1, nou}  that 
\begin{subequations}\label{KKK}
\begin{align}
<e_{z,w},e_{z,w}> & = (1-w\bar{w})^{-2k}\exp(\mc{F}),\\
2\mc{F} &= \frac{2z{\bar{z}}+z^2\bar{w}+\bar{z}^2w}{1-w\bar{w}}
=2|\eta|^2-\bar{w}\eta^2-w\bar{\eta}^2.
\end{align}
\end{subequations}

From \eqref{csv} and \eqref{KKK}, we get for \eqref{etild} on
$\mc{D}^J_1$ the relation
\begin{equation}\label{legPH}
\underline{e}_{\eta,w}=\exp[\frac{1}{4}(w\overline{\eta}^2-cc)]\tilde{e}_{z,w},\quad z=FC(\eta)=\eta-w\bar{\eta}.
\end{equation}

The general scheme \cite{sbcag,sbl} associates to elements of the
Lie algebra $\got{g}$  first order holomorphic differential operators
with polynomial coefficients $X\in\got{g}\rightarrow\db{X}$ as in \eqref{VBC}. The
calculation in \cite{jac1},  based on 
\eqref{BCH}, gives:  
\begin{lemma}\label{mixt}The differential action of the generators
 of
the Jacobi algebra {\em (\ref{baza})} is given by the formulas:
\begin{subequations}\label{summa}
\begin{eqnarray}
& & \mb{a}=\frac{\pa}{\pa z};~\mb{a}^{\dagger}=z+w\frac{\pa}{\pa z} ,
~z,w\in\C, ~|w|<1; \\
 & & \db{K}_-=\frac{\pa}{\pa w};~\db{K}_0=k+\frac{1}{2}z\frac{\pa}{\pa z}+
w\frac{\pa}{\pa w};\\
& & \db{K}_+=\frac{1}{2}z^2+2kw +zw\frac{\pa}{\pa z}+w^2\frac{\pa}{\pa
w} . 
\end{eqnarray}
\end{subequations}
\end{lemma}
Applying Proposition \ref{propvechi} to the representation given by
Lemma \ref{mixt}, we have obtained  in \cite{FC} Berezin's equations
of motion: 
\begin{Proposition}\label{EQLIN1}
The equations of motion on the Siegel-Jacobi disk $\mathcal{D}^J_1$
 generated by
the linear Hamiltonian {\em(\ref{guru})}
 are:
\begin{subequations}\label{qqqN}
\begin{eqnarray}
\ii\dot{z} & = & \epsilon_a+{\bar{\epsilon}}_a w+(\frac{\epsilon_0}{2}
+\epsilon_+  w )z,~z,w\in\C, ~|w|<1, \label{guru2}\\
\ii\dot{w} & = & \epsilon_- + \epsilon_0w+
\epsilon_+w^2 .\label{guru1}
\end{eqnarray}
\end{subequations}
 For the $\eta$ defined in the $FC^{-1}$ transform \emph{(\ref{csv})}, the
 system of first order differential  equations \emph{(\ref{qqqN})} becomes  the
system of separate equations 
\begin{subequations}\label{qqqN1}
\begin{eqnarray}
i\dot{\eta} & = & \epsilon_a
+\epsilon_-\bar{\eta}+\frac{\epsilon_0}{2}\eta,~\eta \in\C,\label{guru22}\\
i\dot{w} & = & \epsilon_- + \epsilon_0w+
\epsilon_+w^2,~ w\in\C,|w|<1, .\label{guru12}
\end{eqnarray}
\end{subequations}
With the change of function $w=XY^{-1}$, the Riccati equation
\eqref{guru1} became the linear hamiltonian system
\begin{equation}\label{GURRU11}
\left(\begin{array}{c}\dot{X} \\ \dot{Y}\end{array}\right)=h_c
\left(\begin{array}{c}X \\ Y\end{array}\right),\quad h_c=\ii 
\left(\begin{array}{cc} -\frac{\epsilon_0}{2} & -\epsilon_-\\
    \epsilon_+  & 
    \frac{\epsilon_0}{2} \end{array}\right) \in \got{sp}(1,\R)_{\C}. 
\end{equation}
If in \eqref{guru22} we make the change of variables $\eta=\xi-\ii
\zeta$, then we  get the system of linear differential equations in real
($\epsilon_a=b+\ii a$)
\begin{equation}\label{GURRU12}
\dot{Z}=h_rZ+F, \quad Z=
\left(\begin{array}{c}\zeta \\ \eta\end{array}\right), ~ F=\left(\begin{array}{c} a \\
    b \end{array}\right),   ~ h_r= 
\left(\begin{array}{cc} n & m-p\\
    m+p  & 
    -n  \end{array}\right) \in \got{sp}(1,\R). 
\end{equation}
\end{Proposition}
We also reproduce the results concerning the phases obtained in \cite{FC}:
\begin{Proposition}\label{prFAZ}
The Berry phase on the Siegel-Jacobi disk $\mc{D}^J_1$ expressed in
the variables $(\eta,w)$, $z=\eta-w\bar{\eta}$, reads
\begin{equation}\label{FFV}
\frac{2}{\ii}d\phi_B=
(\frac{2k\bar{w}}{1-w\bar{w}}-\frac{\bar{\eta}^2}{2}) dw
+(\bar{\eta} + \bar{w}\eta)d \eta - cc . 
\end{equation}
The energy  function \eqref{clH} attached to the Hamiltonian
\eqref{guru}  in  the
coherent state representation \eqref{csu} can be
written as  $\mc{H}=
\mc{H}_{\eta}+\mc{H}_{w}$, where 
\begin{subequations}\label{realHH}
\begin{align}
~\mc{H}_{\eta}  ~ & =  ~  \bar{\epsilon}_a\eta+\epsilon_a\bar{\eta} + 
\frac{1}{2}(\epsilon_+\eta^2+\epsilon_-\bar{\eta}^2+\epsilon_0\eta\bar{\eta}),\label{realHH1} \\
~ \mc{H}_{w} ~ & =  ~  k\epsilon_0+ \frac{2k}{1-w\bar{w}}(\epsilon_+w+\epsilon_-\bar{w}+\epsilon_0w\bar{w}). \label{realHH2}
\end{align}
\end{subequations}
\end{Proposition}
\subsection{Equations of motion on the Siegel-Jacobi ball
  $\mc{D}^J_n$}\label{lasst} 
Following \cite{nou}, we  consider  a  Hamiltonian linear in the generators of
the group $G^J_n$ 
\begin{equation}\label{HACA}
\mb{H}= \epsilon_i\mb{a}_i+\overline{\epsilon}_i\mb{a}_i^{\dagger} +  
\epsilon^0_{ij}\mb{K}^0_{ij}+
\epsilon^-_{ij}\mb{K}^-_{ij}+\epsilon^+_{ij}\mb{K}^+_{ij}. 
\end{equation}
 The hermiticity condition imposes to the matrices of  coefficients $\epsilon_{0,\pm}=(\epsilon^{0,\pm})_{i,j=1,\dots,n} $  the restrictions
\begin{equation}\label{CONDI}
\epsilon_0^{\dagger}=\epsilon_0; ~\epsilon_-=\epsilon_-^t;~
\epsilon_+=\epsilon_+^t; ~ \epsilon_+^{\dagger}=\epsilon_-.
\end{equation}
It is useful to introduce   the matrices $m,n,p,q\in\text{M}(n,\R)$
such that 
\begin{equation}\label{epsmn}
\epsilon_-=m+\ii n, ~
  \epsilon_0^t/2=p+\ii q; p^t=p; m^t= m; n^t =n; q^t=-q.
\end{equation} 
We consider 
 a matrix Riccati equation \eqref{RICC} on the manifold $M$ and a linear differential
equation \eqref{LINZ}  in $z\in\C^n$
  \begin{subequations}\label{TOTAL}
\begin{align}
\dot{W} & =AW+WD+B+WCW, ~A,B,C,D\in M(n,\C); \label{RICC}\\
\dot{z} & = M+Nz; ~M= E+WF; ~ N= A+WC, ~E,F\in C^n. \label{LINZ}
\end{align}
\end{subequations}
Firstly, we recall  how {\bf to solve the matrix Riccati
equation} \eqref{RICC}  {\bf by linearization}.
{\it If we proceed to the homogenous coordinates} $W=XY^{-1}$, $X,Y\in
M(n,\C)$,  {\it a linear
system of ordinary differential equations is attached to the matrix
Riccati equation} \eqref{RICC} (cf. \cite{levin}, see also \cite{sbl})
\begin{equation}\label{RICClin}
\left( \begin{array}{c}\dot{X}\\\dot{Y}\end{array}\right)=h
\left( \begin{array}{c}{X}\\{Y}\end{array}\right), ~
h=\left(\begin{array}{cc} A & B\\ -C & -D \end{array}\right) . 
\end{equation}
{\it Every solution of} (\ref{RICClin}) {\it is a solution of} (\ref{RICC}),
{\it whenever} $\det (Y)\not=0$.

\begin{Proposition} \label{POYT}
The classical motion and quantum evolution generated by  the linear
hermitian  Hamiltonian \eqref{HACA}, \eqref{CONDI}  are described by
first order differential equations:\\
a) On $\mc{D}^J_n$,  $(z,W)\in \C^n\times\mc{D}_n$ verifies
\eqref{TOTAL}, 
with coefficients
 \begin{subequations}\label{hip}
\begin{align}
A_c & =  -\frac{\ii}{2}\epsilon_0^t, ~ B_c=-\ii\epsilon_-, ~C_c=-\ii\epsilon_+, ~
D_c= A_c^t ; \label{hip2}\\
E_c & =-\ii\epsilon, ~F_c=-\ii\bar{\epsilon}\label{hip1}.
\end{align}
\end{subequations}
b) Explicitly,  the differential equations  for $(W,z)\in\mc{D}^J_n$ are
\begin{subequations}\label{Mhip}
\begin{eqnarray}
\ii \dot{W}  & = &\epsilon_-  +
(W\epsilon_0)^s +W\epsilon_+W, ~ W\in\mc{D}_n,\label{Mhip2}\\
\ii \dot{z} & = & \epsilon +
W\overline{\epsilon}+  \frac{1}{2}\epsilon_0^tz+W\epsilon_+z,
~z\in\C^n , \label{Mhip1}
\end{eqnarray}
\end{subequations}
c) Under the $FC$ transform, $z=\eta-W\bar{\eta}$,   
the differential equations 
 in the variables $\eta\in\C^n$, $W\in\mc{D}_n$ become independent:  $W$
verifies \eqref{RICC}  with coefficients \eqref{hip2} and $\eta $ verifies  
\begin{equation}\label{hipPRT1}
\ii \dot{\eta} =  \epsilon +
\epsilon_-\bar{\eta} + \frac{1}{2}\epsilon_0^t\eta,~\eta\in\C^n.
\end{equation}
d) The linear system of differential equations  \eqref{RICClin} attached to the matrix
Riccati equation \eqref{Mhip2} is
 \begin{equation}\label{Rlin}
\left( \begin{array}{c}\dot{X}\\\dot{Y}\end{array}\!\right)=h_c
\left(\! \begin{array}{c}{X}\\{Y}\end{array}\right), ~
h_c=\left(\begin{array}{cc} -\ii(\frac{\epsilon_0}{2})^t &
    -\ii\epsilon_-\\ \ii\epsilon_+& 
\ii\frac{\epsilon_0}{2}\end{array}\right)\in \got{sp}(n,\R)_{\C} , ~W=X/Y\in\mc{D}_n. 
\end{equation}
e) In \eqref{hipPRT1} we
introduce $\eta=\xi - \ii \zeta$, $\xi,\zeta\in\R^n$ and  we put 
$\epsilon= b+\ii a$, where $a,b\in\R^n$.  The first order complex differential equation equation
\eqref{hipPRT1} is equivalent with a system of first order real 
differential equations with real coefficients, which we write as  
\begin{equation}\label{LINe} \dot{Z}=h_rZ + F, ~ Z = 
\left( \begin{array}{c}{\xi} \\
    {\zeta}\end{array}\right), ~ 
 F = \left(\begin{array}{c}a \\ b \end{array}\right),
h_r\!=\! \left(\begin{array}{cc} \!n+q
    & m-p\\  m+ p  & -n+q \!\end{array}\!\right)\! \in\got{sp}(n,\R).
\end{equation}
\end{Proposition}
\vspace{3ex}

\subsection*{Acknowledgments}
This research  was conducted in  the  framework of the 
ANCS project  program PN 09 37 01 02/2009  and     the UEFISCDI - Romania 
 program PN-II Contract No. 55/05.10.2011.

\today
\end{document}